\documentclass[12pt]{article}
\usepackage[margin=1in]{geometry}
\usepackage{graphicx}
\usepackage{hyperref}
\usepackage{tcolorbox}
\usepackage{enumitem}
\usepackage[authoryear]{natbib}
\usepackage{fullpage}
\usepackage{xcolor}
\usepackage{amsmath}
\usepackage{amssymb}
\usepackage{bold-extra}
\usepackage{amsthm}
\usepackage{mathtools}
\usepackage{setspace}
\usepackage{appendix}
\usepackage{float}
\usepackage{authblk}

\theoremstyle{remark}

\newtcolorbox{convergencebox}{
    colframe=black!75!white,
    colback=gray!10,
    sharp corners,
    boxrule=0.8mm,
    title={\textbf{Convergence Analysis}},
    colbacktitle=black!70,
    fonttitle=\bfseries
}

\newtcolorbox{algorithmbox}[1]{colframe=black!75!white,
		colback=gray!10, sharp corners=south, boxrule=0.8mm,
		title={\textbf{#1}}, colbacktitle=black!70, fonttitle=\bfseries}

\newcommand{\Ab}{\mathbb{A}}
\newcommand{\Eb}{\mathbb{E}}
\newcommand{\Hb}{\mathbb{H}}
\newcommand{\Kb}{\mathbb{K}}

\newcommand{\Pb}{\mathbb{P}}

\newcommand{\Xb}{\mathbb{X}}
\newcommand{\Yb}{\mathbb{Y}}

\newcommand{\Bc}{\mathcal{B}}

\newcommand{\Xc}{\mathcal{X}}
\newcommand{\Hc}{\mathcal{H}}
\newcommand{\Tc}{\mathcal{T}}

\newcommand{\A}{\mathbb{A}}

\newcommand{\K}{\mathbb{K}}

\newcommand{\R}{\mathbb{R}}

\newcommand{\X}{\mathbb{X}}
\newcommand{\Y}{\mathbb{Y}}

\renewcommand{\bar}{\overline}

\newcommand{\defeq}{:=}

\DeclarePairedDelimiter{\set}{\{}{\}}

\title{Computing optimal policies for managing inventories with noisy observations\footnote{This research was partially supported by the U.S. Office of Naval Research (ONR) under
grants N000142412608, N000142412646, and N0001425GI01179.}}

\author[1]{Eugene Feinberg}
\author[2]{Jefferson Huang}
\author[3]{Pavlo Kasyanov}
\author[1]{Thomas O'Neill}

\affil[1]{Stony Brook University}
\affil[2]{Naval Postgraduate School}
\affil[3]{Institute for Applied System Analysis, National Technical University of Ukraine, ``Igor Sikorsky Kyiv Polytechnic Institute''}

\date{}

\begin{document}

\maketitle

\begin{abstract}
 This paper implements the Deep Deterministic Policy Gradient (DDPG) algorithm for computing optimal policies for partially observable single-product periodic review inventory control problems with setup costs and backorders. The decision maker does not know the exact inventory level, but can obtain noise-corrupted observations of them.  The goal is to maximize the  expected total discounted costs incurred over a finite planning horizon. We also investigate the Gaussian version of this problem with normally distributed initial inventories, demands, and observation noise. We show that expected posterior observations of inventory levels, also called mean beliefs,  provide sufficient statistics for the Gaussian problem. Moreover, they can be represented in the form of a Markov Decision Processes for an inventory control system with time-dependent holding costs and demands.  Thus, for a Gaussian problem, the there exist $(s_t,S_t)$-optimal policies based on mean beliefs, and this fact explains the structure of the approximately optimal policies computed by DDPG. For the Gaussian case, we also numerically compare the performance of policies derived from DDPG to optimal policies for discretized versions of the original continuous problem.  
\end{abstract}

\section{Introduction} \label{s2} Decision making with incomplete information is an important topic that is challenging from both computational and theoretical points of view.  One of the most popular models, a Partially Observable Markov Decision process, is computationally challenging even in the case of finite state, action, and observation sets; see  \citep{PapTsi}, and additional difficulties arise when one or more of these sets are infinite. The main solution approach is to reduce the problem to a so-called belief Markov Decision Process (MDP) whose states are posterior probability distributions over states that depend on observations. For a finite-state problem with $n$ states, the state space of the belief MDP is the $(n-1)$-dimensional probability simplex. For infinite-state problems, the state space of the belief MDP is also infinite-dimensional. This creates additional theoretical and computational challenges. Since the problem is difficult, the natural approach is to use reinforcement learning (RL); see \cite{SuttonBarto}.

In this paper, we consider a partially observable version of the classic single-product periodic-review inventory problem with setup costs and backorders. Here, the Decision Maker (DM) cannot observe the actual inventory level. Instead, the DM only observes the inventory level after it has been corrupted by additive noise. The goal is to minimize the expected discounted costs over a finite planning horizon. We consider the Deep Deterministic Policy Gradient (DDPG) algorithm proposed by \citep{silver2014deterministic} specialized to this partially observable problem, with states consisting of action-observation histories. The convergence of this algorithm was observed empirically, and the policies that it produced were similar to $(s,S)$ policies.

We also investigated a Gaussian version of the problem when the initial inventory control levels, demands, and observation noises have normal distributions. If the objective function where quadratic, then this would be a classic partially observable linear-quadratic control problem.  However, the objective function for the inventory control problem is not quadratic. For the Gaussian version, the belief probability distributions are normal distributions, which means they are parametrized by their mean and variance.  So, instead of dealing with measures on the real line, the state of the belief MDP can be represented by a two dimensional vector
$(\bar{x},\sigma^2)$ representing the mean and variance of beliefs.  For the inventory control problem, this MDP has remarkable properties.  If we consider a sequence of belief states $\{(x_t,\sigma_t^2)\}_{t=0}^\infty$ at times $t=0,1,\ldots,$ then \{$\sigma_t\}_{t=0}^\infty$ is a deterministic sequence depending only on the parameters of the original problem and the initial state distribution. This means we can instead consider an MDP where the state $\bar{x}$ corresponds to the mean of the current posterior state distribution. In addition, this MDP corresponds to a single-product periodic review inventory control problem with setup costs, backorders, and where the holding costs and demands are time-dependent.  For such problems, there exist optimal $(s_t,S_t)$ policies.  We also applied the DDPG algorithm to this problem, and produced optimal policies that are functions of expectations of beliefs corresponding to observed histories.  The one-dimensional problem was also solved by state discretization; see e.g., \cite{saldi2018finite}.

The discovered reduction of the belief MDP to a single-product nonstationary inventory problem can be also applied to infinite-horizon average costs per unit time.  For average-costs per unit time, solving even simple inventory control problems with incomplete information is a nontrivial problem \cite{} because it is difficult to apply the sufficient conditions for average-cost optimality developed by \cite[Assumptions (B), (${\rm \underline{B}}$)]{FKZ} to infinite-dimensional problems. However, we conjecture that the described periodic-review inventory control problems for average beliefs also have optimal $(s_t,S_t)$ policies that correspond to an optimal policy for the original problem.

The problem formulation is provided in Section~\ref{s2}.  Section~\ref{s3} provides the DDPG algorithm description for the original problem where states are observation histories.
Sections~\ref{s4} and \ref{s5} provide descriptions of the DDPG and discretization algorithms for the Gaussian version of the problem.  Section~\ref{s6} provides the descriptions of numerical inputs and computational results.
The appendices describe mathematical issues relevant to this paper. \ref{secA1} describes a reduction of discrete-time  stochastic filtering problems to POMDPs.  The inventory control problem with noisy observations studied in this paper is an example of such a problem. \ref{secA2} describes the construction of a belief MDP. \ref{secA3} describes the reduction of the belief MDP for an inventory control problem with noisy observations studied in this paper to a single-product inventory problem with time-dependent demands and holing costs, and without noise. This reduction explains the results produced by the DDPG algorithm. In a broader sense, this paper illustrates the recently understood links in \citep{FIKK} between POMDPs and statistical filtering by solving a partially-observable version of a classic inventory control problem.

\section{Periodic-Review Inventory Control with Noisy Observations} \label{s2}
Consider a discrete-time setting with time steps $t = 0, 1, 2, \dots$. At the start of each step $t$,  the current inventory level $x_t$ takes values in $\mathbb{R}:= (-\infty,+\infty).$ The decision maker (DM) does not know the exact level of inventory because of observation/recording errors. At each step $t,$ the DM selects an order quantity $a_t\in\mathbb{R}_{\ge 0}:=[0.+\infty).$ The order is received instantly and added to the inventory.
	
Following the placement of the order, the random demand $D_t$ is realized. The sequence $\{D_t\}_{t \geq 1}$ consists of independent and identically distributed nonnegative random variables with the distribution function $F_D.$ The inventory level at the start of the next step is
\[
x_{t+1} = L( x_t + a_t - D_t),\qquad t=0,1,\ldots,
\]
where, for $w\in\R,$ $L(w):=w$ for problems with back orders, and $L(w):=w^+:=\max\{w,0\}$ for problems with lost sales. At the initial time step $t=0,$ the initial state $x_0$ has a known probability distribution $p_0.$ However, the inventory levels $x_t$ are measured with errors, and the DM knows only the values
\[ y_{t}=x_{t}+\eta_{t}, \qquad t=0,1,\ldots,\]
where $\eta_0,\eta_1,\ldots$ are iid random variables with the distribution function $F_\eta.$

At the end of period $t$, the system incurs a cost consisting of two components: a holding cost and an ordering cost. If at a state $x$ an action $a$ is selected, then the ordering cost is

\[ \hat{c}(a)=K\cdot\mathbf{1}_{a>0}+\tilde{c}a, \]
where
$K$ is a constant representing the fixed ordering cost, while $\tilde{c}$ represents the per-unit ordering cost. 
The holding cost depends on the actual inventory level after orders and demands are realized, and is represented by a nonnegative convex function $h$.
If a decision $a$ is chosen at a state $x,$ then the expected one-stage cost is

\[ c(x,a)=\hat{c}(a)+\mathbb{E}[h( L(x+a-D))], \]
 where $D\sim F_D.$

We denote by $\X$   the set of possible inventory levels, and by $\Y$ the set of possible observations, where $\X:=\R$ and $\Y:=\R.$  Let $\Hb_t= (\Y\times\A)^t\times \Y$ be the set of all possible histories that can be observed by the DM before an action is chosen at step $t,$.  At each step $t=0,1,\ldots,$ the DM chooses an action $a_t$ according to decision a rule $\pi_t(da_t|p_0,y_0,a_0,y_1,\ldots,a_{t-1},y_t),$  where a decision rule $\pi_t$ is a regular transition probability from $\Hb_t$ to $\A.$ A nonradomized decision rule at step $t$ is a Borel mapping $\varphi_t:\Hb_t\to \A.$  A policy $\pi$ is a sequence of decision rules $\pi=(\pi_0,\pi_1,\ldots).$ The set of all policies is denoted by $\Pi.$   A nonrandomized policy $\varphi$ is a sequence of decision rules $\varphi=(\varphi_0,\varphi_1,\ldots).$  Every policy $\pi$ defines a probability measure $P^\pi$ on the set of trajectories $(\X\times\Y\times\A)^\infty,$ and we define by $\Eb^\pi$ expectations with respect to this measure. The set of all policies is denoted by $\Pi.$

For a finite horizon $T=1,2,\ldots$ and for an infinite horizon $T=\infty,$ the expected total discounted cost is
\begin{equation}\label{eqdisccost}
v_{\alpha,T}^\pi:=\Eb^\pi\sum_{t=0}^{T-1}\alpha^t c(x_t,a_t),
\end{equation}
where $\alpha\in [0,1)$ is the discount factor. Let us consider the optimal objective values $v_{\alpha,T}=\inf_{\pi\in\Pi}v_{\alpha,T}^\pi.$  A policy $\pi$ is  optimal (or discount optimal) for the horizon $T$ if $v_{\alpha,T}^\pi=v_{\alpha,T}.$ For infinite-horizon problems we usually use the notations $v_{\alpha}^\pi:=v_{\alpha,\infty}^\pi$ and $v_{\alpha}:=v_{\alpha,\infty}.$ For computational purposes, it is useful to use $Q$-functions defined as ${\bf Q}_{\alpha,t}(x,a)=\mathsf{E}[c(x_0,a_0)+v_{\alpha,  {t-1}}(x_1)|x_0=x,a_0=a].$

According to \cite[Corollary 7.6]{FIKK}, if the observation noises $\eta_t$ are continuous random variables, then for each horizon $T=1,2,\ldots$ or $T=\infty$ there exists an optimal policy.  As explained in \cite{FIKK}, this is true because the original problem can be viewed as a partially observable Markov decision process (POMDP), and this POMDP can be reduced to a belief-MDP with weakly continuous transition probabilities and $\K$-inf-compact one-step costs. States in the belief-MDP correspond to probability distributions of states in the original problem.  For this belief-MDP there exist optimal Markov policies for finite horizon problems, and there exists an optimal deterministic policy for the infinite-horizon problem.  These optimal policies for the belief-MDP define optimal policies for the original POMDP; see e.g., \cite{HL}.  In addition, it is possible to write optimality equations for the belief-MDPs.  These optimality equations define Markov optimal policies for finite-horizon belief-MDPs and deterministic optimal policies for an infinite-horizon belief-MDPs.  In addition, value iterations applied to finite-horizon values converge to infinite-horizon values.  However, this approach is computationally inefficient; therefore, we consider other approaches based on reinforcement learning, and use the theory of POMDPs to explain the results.

We also consider a Gaussian model with normal initial state distribution $p_0,$ demand distribution $F_D,$ and noise distribution $F_\eta.$  In particular, $F_D\sim \mathcal{N}(\bar{D},\sigma_D^2)$ and  $F_\eta\sim \mathcal{N}(\bar{D},\sigma_\eta^2).$  Of course, demand should be nonnegative, but from the practical point of view this can be achieved if $\bar{D}>0$ and $\sigma_D$ is small relative to $\bar{D}$.

\section{Applying Reinforcement Learning}\label{s3}

The partially observable inventory control problem described in the previous section is computationally intractable. Even though it is possible to reduce it to a belief-MDP, the state space is uncountably infinite. Moreover, while computational methods exist  for POMDPs with finite state and action sets (see \cite{smallwood1973optimal}), the POMDP corresponding to the partially observable inventory control problem has uncountably infinite state and action sets.

In this paper, we adapt the DDPG algorithm (\cite{silver2014deterministic}) to the partially observable inventory control problem described in Section~\ref{s2}. DDPG is a policy gradient algorithm (see \cite{sutton2000policy}) that is designed specifically for MDPs with continuous state and action sets. For the inventory control problem, we consider policies learned via DDPG with neural network function approximation; specifically, we use multi-layer perceptrons with two hidden layers to represent both the (deterministic) policy (the ``Actor'') as well as the $Q$-function (the ``Critic''). We evaluated DDPG in both Gaussian and non-Gaussian problem contexts. For the non-Gaussian case, we apply DDPG to an MDP associated with the original partially observable problem, where the states of the MDP correspond to observed histories. For the Gaussian case, the associated belief-MDP can be formulated so that the state is a scalar corresponding to the mean of the current posterior state distribution; see \ref{secA3}. Therefore, for the Gaussian case we compare the performance of policies learned by DDPG on this belief-MDP to that of policies learned by applying DDPG to the MDP with histories as states. We find that both approaches produce policies that perform comparably to policies obtained by discretizing the belief-MDP and solving the resulting finite MDP via value iterations; see Section~\ref{s6}. However, applying DDPG to the belief-MDP requires much less computational effort than both applying DDPG to the MDP with histories and solving the discretized MDP. Moreover, the state dimension in the MDP with histories increases as the planning horizon increases, while in the Gaussian case the state of the belief MDP remains a scalar regardless of the planning horizon length.

\subsection{Pseudocode} \label{pseudocode}

In this sub-section, we provide pseudocode that describes the modifications we made to DDPG to adapt it to the partially observable inventory control problem that we consider. This pseudocode uses the following definitions and subroutines:

\paragraph{Definitions}
\begin{itemize}
    \item replay buffer $\mathcal{B}$ \textemdash stores vectors $(hist_{t+1}, c, done)$ for training
    \item learning rate \textemdash size of update to network parameters during each optimization step
    \item batch\_size \textemdash sample size for calculating expectations during training
    \item max\_action \textemdash maximum action allowed during simulation (reduces time spent learning states that will likely never be reached by the optimal policy)
    \item update rate \textemdash rate at which target network parameters adapt to non-target network parameters
    \item discretized\_xvalues \textemdash discrete set of points at which to evaluate belief distribution
    \item min\_zvalue \textemdash minimum belief probability density value to save (improves efficiency)
    \item $\delta_c$ is the random variable equal to $c$ with probability 1.
    \item $hist_t$ \textemdash history at time $t$ stored as a vector of length $2\cdot\text{len\_episode}+2$ as follows: $[t,x_0,a_0,y_1,a_1,y_2,\dots,a_t,y_{t+1},0,0,\dots0]$
    \item $cti$ \textemdash critic training iterations
    \item $ati$ \textemdash actor training iterations
    \item $lr$ \textemdash neural network learning rate
    \item $\beta_1$ \textemdash coefficient used for computing running averages of gradient \citep{pytorch_adam}
    \item $\beta_2$ \textemdash coefficient used for computing running averages of square of gradient \citep{pytorch_adam}
    \item steps\_done \textemdash total steps completed over all episodes\
    \item $\epsilon$ \textemdash probability to explore rather than choose optimal action
    \item $\epsilon_s$ \textemdash starting value of $\epsilon$
    \item $\epsilon_e$ \textemdash ending value of $\epsilon$
    \item eps\_decay \textemdash $\epsilon$ decay factor
    \item exploration\_noise \textemdash variance of normal distribution added to chosen action when the agent explores
    \item $\mathbb{X}^\mathcal{D}_t$ \textemdash quantized set of inventory levels with belief probability exceeding min\_zvalue at a given time $t$

\end{itemize}

\paragraph{Subroutines}
\begin{itemize}
\item NextBelief($z_t,a_t,y_{t+1}$) \textemdash Takes the most recent belief, observation, and action
\item Train($\mathcal{B}$) \textemdash
Trains the actor and critic networks based on data from the replay buffer
\item Optimize($L$,$\Omega_{i-1}$) \textemdash We used the PyTorch implementation \cite{pytorch_adam} of the Adam \citep{kingma2014adam} optimization algorithm during the neural network training process
\end{itemize}

We provide the pseudocode for the subroutines NextBelief($z_t,a_t,y_{t+1}$) and Train($\mathcal{B}$) after the pseudocode describing our modifications to DDPG, which follows:

\begin{algorithmbox}{DDPG for Inventory Control}\footnotesize

    Input: actor parameters $\theta$, target actor parameters $\theta'$, critic parameters $\phi$, target critic parameters $\phi'$ \hfill $\triangleright$ Initial network parameters (arbitrary) \\
    $\mathcal{B}\gets\emptyset$ \hfill $\triangleright$ Initialize empty replay buffer \\
    steps\_done$\gets0$ \hfill $\triangleright$ Initialize step counter \\

    \vspace{6pt}
    Repeat until policy $\mu_\theta$ converges:\\
    \hspace*{1.5em} $x_0 \sim p_0(x)$ \hfill $\triangleright$ Sample initial state\\
    \hspace*{1.5em} Sample $\eta_0 \sim F_\eta$; compute $y_{0} = x_{0} + \eta_0$ \hfill $\triangleright$ Observe noisy initial state\\
    \hspace*{1.5em} $hist_0\gets[0,y_0,\dots,0,\dots,0]$ \hfill $\triangleright$ Initialize history vector\\
    \hspace*{1.5em} $\mathbb{X}^{\mathcal{D}} \gets \text{discretized\_xvalues}$, $z(x) \gets p_0(x), x\in \mathbb{X}^{\mathcal{D}}$ \hfill $\triangleright$ Initialize belief\\
    \hspace*{1.5em} \textbf{For} each timestep $t = 0, 1, \dots, \text{len\_episode}-1$\\
    \hspace*{3em} $\epsilon = \epsilon_e + (\epsilon_{s}-\epsilon_e) * \exp(-\frac{\text{steps\_done}}{\text{eps\_decay}})$\\
    \hspace*{3em} $a_t \sim (1-\epsilon)\delta_{\mu_{\theta}(hist)} + \epsilon \delta_{\mathcal{N}(0,\text{exploration\_noise})}$ \hfill $\triangleright$ Choose action with exploration noise\\
    \hspace*{3em} Clip $a_t$ to $[0, \text{max\_action}]$ \hfill $\triangleright$ Enforce action bounds\\
    \hspace*{3em} Sample $D_t \sim F_D$; compute $x_{t+1} = x_t + a_t - D_t$ \hfill $\triangleright$ Inventory update\\
    \hspace*{3em} Sample $\eta_t \sim F_\eta$; compute $y_{t+1} = x_{t+1} + \eta_{t+1}$ \hfill $\triangleright$ Observe noisy state\\
    \hspace*{3em} $hist_{t+1} \gets [t+1,y_0,a_0,y_1,a_1,y_2,\dots,a_t,y_{t+1},0,0,\dots0]$ \hfill $\triangleright$ Update history\\
    \hspace*{3em} $z_{t+1}, \mathbb{X}^{\mathcal{D}}_{t+1} \gets \text{NextBelief}(z_t, a_t, y_{t+1})$ \hfill $\triangleright$ Belief update\\
    \hspace*{3em} \textbf{If} terminal state ($t = \text{len\_episode}-1$): $done \gets 1$, else $done \gets 0$ \hfill $\triangleright$ Set terminal flag\\
    \hspace*{3em} Compute $c(hist_t,a_t) = \sum_{x \in \mathbb{X}^{\mathcal{D}}} c(x,a_t) \cdot z(x)$ \hfill $\triangleright$ Expected cost\\
    \hspace*{3em} Store $(hist_{t+1}, c(hist_t,a_t))$ in $\mathcal{B}$ \hfill $\triangleright$ Save transition to buffer\\
    \hspace*{3em} Train using $\text{Train}(\text{batch\_size})$ \hfill $\triangleright$ Update actor and critic\\
    \hspace*{3em} steps\_done$\gets$steps\_done$+1$

\end{algorithmbox}

\vspace{0.5cm}

\begin{algorithmbox}{NextBelief($z_t,a_t,y_{t+1}$) Function}\footnotesize
    \textbf{For} $x_2$ in discretized\_xvalues\\
    \hspace*{1.5em} $z^\text{unnormalized}_{t+1}(x_2)\gets f_{\eta}(y - x_2)\sum_{x_1\in \mathbb{X}^\mathcal{D}} f_D(x_1 + a - x_2) z(x_1)$ \hfill $\triangleright$ Bayesian update (unnormalized)\\

    $z_{t+1}^{area}\gets\frac{1}{2}\sum_{x_i\in\text{discretized\_xvalues}} (z_{t+1}(x_{i+1}) - z_{t+1}(x_i))(x_{i+1} - x_i)$ \hfill $\triangleright$ Approximate area under belief\\

    $z_{t+1}(x)\gets z_{t+1}^\text{unnormalized}(x)/z_{t+1}^{area}$ \hfill $\triangleright$ Normalize untrimmed belief\\

    $\mathbb{X}^{\mathcal{D}}_{t+1} \gets [x \in \text{discretized\_xvalues} \mid z_{t+1}(x) > \text{min\_zvalue}]$ \hfill $\triangleright$ Trim low-probability states\\

    $z_{t+1}^\text{trimmed\_area} \gets \frac{1}{2}\sum_{x^{(i)} \in \mathbb{X}^{\mathcal{D}}} (z_{t+1}(x^{(i+1)}) - z_{t+1}(x^{(i)}))(x^{(i+1)} - x^{(i)})$ \hfill $\triangleright$ Recalculate area after trimming\\

    $\forall x \in \mathbb{X}^\mathcal{D}_{t+1}: z_{t+1}(x) \gets z_{t+1}(x)/z_{t+1}^\text{trimmed\_area}$ \hfill $\triangleright$ Normalize trimmed belief\\

    \textbf{return} $z_{t+1}, \mathbb{X}^\mathcal{D}_{t+1}$ \hfill $\triangleright$ Output updated belief and support\\
\end{algorithmbox}

\begin{algorithmbox}{Train($\mathcal{B}$) Function}\footnotesize
    If $|\mathcal{B}| < \text{batch\_size}$, exit function. \hfill $\triangleright$ Wait for enough samples\\

    Sample mini-batch of size $\text{batch\_size}$ of $(hist_{t+1}, c)$ from $\mathcal{B}$ \hfill $\triangleright$ Draw training data\\
    ($hist_{t+1}$ contains $hist_t$ and $a_t$, so $hist_{t+1}$ defines $hist_t$, $a_t$, and $c(hist_t,a_t)$)\\

    Perform $cti$ times:\\
    \hspace*{1.5em}$L_{\text{critic}} \gets \mathbb{E}_{\mathcal{B}} \left[(\textbf{Q}_{\phi}(hist_t, a_t) - [c + (1 - done) \cdot \gamma \cdot \textbf{Q}_{\phi'}(hist_{t+1}, \mu_{\theta'}(hist_{t+1}))])^2\right]$ \hfill $\triangleright$ Critic loss (Bellman error)\\
    \hspace*{1.5em}$\phi \gets \text{Optimize}(L_{\text{critic}}, \phi)$ \hfill $\triangleright$ Update critic parameters\\

    Perform $ati$ times:\\
    \hspace*{1.5em}$L_{\text{actor}} \gets \mathbb{E}_{\mathcal{B}}[\textbf{Q}_{\phi}(hist_t, \mu_{\theta}(hist_t))]$ \hfill $\triangleright$ Actor objective (maximize value)\\
    \hspace*{1.5em}$\theta \gets \text{Optimize}(L_{\text{actor}}, \theta)$ \hfill $\triangleright$ Update actor parameters\\

    Soft-update target network parameters: \\
    $\theta' \leftarrow \tau \theta + (1 - \tau) \theta'$, $\phi' \leftarrow \tau \phi + (1 - \tau) \phi'$ \hfill $\triangleright$ Target smoothing\\
\end{algorithmbox}

We remark that our adaptation of DDPG to partially observable inventory control can be modified for use in the completely observable setting by replacing $hist$ with the true current inventory level.

\section{The Gaussian Case}\label{s4}

When the demand and noise distributions as well as the noise distribution are  Gaussian, the computation of posterior belief states can be significantly simplified. This is because in this setting the posterior belief states are normal distributions. This means the belief MDP can be presented as an MDP with the states being two-dimensional vectors representing the mean and variance of the posterior distribution of the current inventory level; the details of the transformation from the original partially-observed MDP to the belief MDP are provided in \ref{secA3}.

\subsection{Applying DDPG}

In the Gaussian case, we can apply the DDPG algorithm as before, where the only change is that we replace $hist$ with the current mean and standard deviation ($\hat{x}_t,\hat{\sigma}^2_t$). Compared to tracking the entire history, the dimension of the state is reduced from $2\cdot$len\_episode$+2$ to 2. In addition, as explained in \ref{secA3}, this dimension can be reduced to one by excluding variances since the sequence of variances is deterministic and depends neither or actions nor on observations.

\subsection{Discretizing the Gaussian Model}
\label{s5}

We compare the policy computed via DDPG to policies obtained by discretizing the state space and applying value iterations. This is not computationally feasible in the non-Gaussian case, where the posterior state distributions may correspond to arbitrary probability distributions on the real line. In the Gaussian case, we approximate the original state space $\mathbb{X}$ with a discrete state space $\dot{\mathbb{X}}$ that has $m$ evenly spaced values $x_i$ ranging from $\dot{x}_{1}$ to $\dot{x}_{m}$.

For each time step $t = 0, 1, \dots$, the variance $\sigma_t^2$ of the posterior state distribution does not depend on observations. In particular, letting $\sigma_0^2$ denote the variance of the initial (Gaussian) state distribution,
$$
    \sigma_{t}^2 = \begin{cases}
        \frac{\sigma_0^2 \sigma_\eta^2}{\sigma_0^2 + \sigma_\eta^2} & \text{if} \ t = 0 \\
        \frac{(\sigma_{t-1}^2 + \sigma_D^2)\sigma_\eta^2}{\sigma_t^2 + \sigma_D^2 + \sigma_\eta^2} & \text{otherwise}.
    \end{cases}
$$
A derivation of this formula for $\sigma_t^2$ is provided in \ref{secA3}. If the mean of the posterior state distribution on time step $t$ is $\dot{x}_i$, $i \in \{1, \dots, m\}$, then the probability that the next posterior state distribution has a mean of $\dot{x}_j$ is set to the probability that a normally distributed random variable with a mean of $\mu = \frac{\dot{x}_i + \dot{x}_{i+1}}{2} + a - \overline{D}$ and a variance of $\sigma^2 = \hat{\sigma}_t$ falls within the interval $[\dot{x}_j, \dot{x}_{j+1})$ for $j = 1, \dots, m-1$. Specifically, letting $\Phi$ denote the cumulative distribution function of a standard normal random variable, the transition probabilities in the discretized MDP are defined for $i \in \{1, \dots, m\}$, $a \in \mathbb{A}$, and $j \in \{2, \dots, m-1\}$ as follows:
$$
    \dot{p}(\dot{x}_j | \dot{x}_i, a) := \Phi\left(\frac{\dot{x}_{j+1} - (\frac{\dot{x}_i + \dot{x}_{i+1}}{2} + a - \overline{D})}{\sigma_t}\right) - \Phi\left(\frac{\dot{x}_j - (\frac{\dot{x}_i + \dot{x}_{i+1}}{2} + a - \overline{D})}{\sigma_t}\right)
$$
Transitions in the original problem that correspond to transitions outside of the interval $[\dot{x}_1, \dot{x}_m]$ are treated as transitions to either $\dot{x}_1$ or $\dot{x}_m$ by setting
$$
    \dot{p}(\dot{x}_1 | \dot{x}_i, a) := \Phi\left(\frac{\dot{x}_1 - (\frac{\dot{x}_i + \dot{x}_{i+1}}{2} + a - \overline{D})}{\sigma_t}\right), \quad i \in \{1, \dots, m\}, \ a \in \mathbb{A}
$$
and
$$
    \dot{p}(\dot{x}_m | \dot{x}_i, a) := 1 - \Phi\left(\frac{\dot{x}_m - \frac{\dot{x}_i + \dot{x}_{i+1}}{2} + a - \overline{D})}{\sigma}\right), \quad i \in \{1, \dots, m\}, \ a \in \mathbb{A}.
$$
Finally, to define the one-step costs $\dot{c}$, let $X$ have a normal distribution with a mean of $\mu = \frac{\dot{x}_i + \dot{x}_{i+1}}{2} + a - \overline{D}$ and a variance of $\sigma^2 = \hat{\sigma}_t$. Then $\dot{c}(\dot{x}_i, a)$ is defined as the conditional mean of $c(X, a)$ given that $X \in [\dot{x}_i, \dot{x}_{i+1})$.

\section{Numerical Results} \label{s6}   

This section  provides numerical results for non-Gaussian and Gaussian versions of the partially observable inventory control problem.  In the non-Gaussian case, we applied the DDPG algorithm to an MDP with states that correspond to observed histories. In the Gaussian case, we compared three methods: (i) the DDPG algorithm applied to an MDP with states that correspond to observed histories, (ii) solving an MDP obtained via state discretization with value iterations, and (iii) the DDPG algorithm applied to an MDP with states that correspond to  mean beliefs (i.e., the mean of the posterior distribution of the inventory level).

For the non-Gaussian case, the per-period demand $D$ has an exponential distribution with a mean of 1, the observation noise has a standard normal distribution (mean of 0 and variance of 1), and  the initial inventory control level is normally distributed with a mean of 2 and a variance of 4. The cost parameters are the same as in the Gaussian model.  The following table displays the parameters of the Gaussian model and the values of the constants that were used in the implemented algorithms.

\begin{table}[H]\small
    \centering
    \begin{tabular}{|l|l|}
        \hline
        \textbf{Parameter} & \textbf{Value} \\
        \hline
        episode length ($\text{len\_episode}$) & 4 \\
        batch size ($\text{batch\_size}$) & 512 \\
        maximum action ($\text{max\_action}$) & 12 \\
        exploration noise ($\text{exploration\_noise}$) & 8 \\
        update rate ($\tau$) & 0.005 \\
        discretized x-values & $\{m + i\Delta : i = 0, 1, \ldots, (M - m)/\Delta\}$ \\
        discretization lower bound ($m$) & -20 \\
        discretization upper bound ($M$) & 30 \\
        discretization step ($\Delta$) & 0.5 \\
        minimum z-value ($\text{min\_zvalue}$) & 0.0001 \\
        fixed order cost ($K$)& 1 \\
        unit order cost ($\bar{c}$) & 0.1 \\
        initial state distribution ($p_0(x)$) & $\mathcal{N}(2, 4)$ \\
        demand distribution ($F_D$) & $\mathcal{N}(\overline{D}, \sigma_D^2)$ \\
        $\overline{D}$ & 1 \\
        $\sigma_D^2$ & 1 \\
        noise distribution ($F_{\eta}$) & $\mathcal{N}(0, \sigma_{\eta}^2)$ \\
        $\sigma_{\eta}^2$ & 1 \\
        holding cost ($h(x)$) & $\begin{cases}
                x & x>0 \\
                -5x & x\le0
            \end{cases}$ \\

        actor training iterations ($ati$) & 1 \\
        critic training iterations ($cti$) & 3 \\
        actor learning rate ($lr_{actor}$) & $10^{-5}$ \\
        target actor learning rate $lr_{target\_actor}$ & $10^{-5}$ \\
        critic learning rate ($lr_{critic}$) & $10^{-3}$ \\
        target critic learning rate $lr_{target\_critic}$ & $10^{-3}$ \\
        Adam algorithm coefficient ($\beta_1$) & 0.999 \\
        Adam algorithm coefficient ($\beta_2$) & 0.999 \\
        Starting exploration probability ($\epsilon_s$) & 0.9 \\
        Ending exploration probability ($\epsilon_e$) & 0.05 \\
        Exploration decay rate (eps\_decay) & 200 \\

        \hline
    \end{tabular}
    \caption{Model Parameters}
    \label{tab:modelparameters}
\end{table}

The computation time and resulting performance for each algorithm is presented in Table \ref{tbl:results}. The reported value is the average over 3000 episodes of the total discounted cost that was incurred under the given policy; the corresponding standard errors are reported as well. In the Gaussian case, we observe that as the quantization becomes finer (i.e., as dx decreases), the average total discounted cost decreases as expected, with the finest quantization level providing the best estimate for the performance of an optimal policy for the original problem. For this Gaussian case, both DDPG applied to histories as well as DDPG applied to beliefs performed comparably to the policy obtained from the finest quantization level. This suggests that in both cases, DDPG produced a nearly-optimal policy for the original problem. Moreover, DDPG applied to beliefs achieved this level of performance using roughly 20\% less computational time than DDPG applied to histories. This illustrates the benefit of using the structure of the problem (i.e., that the belief distribution can be identified with a scalar corresponding to its mean) to design the MDP to which DDPG is applied to solve the original partially-observed problem.





\begin{table}[H]

    \begin{tabular}{ |c|c|c|c|c| }
        \hline
        \textbf{Problem} & \textbf{Method} & \textbf{Time} & \textbf{Value} & \textbf{Std.\ Error} \\ \hline
        Gaussian & Quantized dx=1.0 & 91.6s & 8.4 & 0.048\\
        \hline

        Gaussian & Quantized dx=0.5 & 687.7s & 8.0 & 0.035\\
        \hline

        Gaussian & Quantized dx=0.3 & 2689.8s & 7.9 & 0.034\\
        \hline

        Gaussian & DDPG with Histories - 15000 episodes & 2093.3s
        & 7.8 & 0.032\\
        \hline
        Gaussian & DDPG with Beliefs - 15000 episodes & 1670.2s
        & 7.9 & 0.033\\
        \hline

        Non-Gaussian & DDPG with Histories - 15000 episodes & 10895.1s
        & 11.1 & 0.033\\
        \hline

    \end{tabular}
    \caption{Computational results for the Gaussian and non-Gaussian cases.}
    \label{tbl:results}
\end{table}

\small

\begin{figure}

\centering
\begin{tabular}{cc}

  \includegraphics[width=65mm]{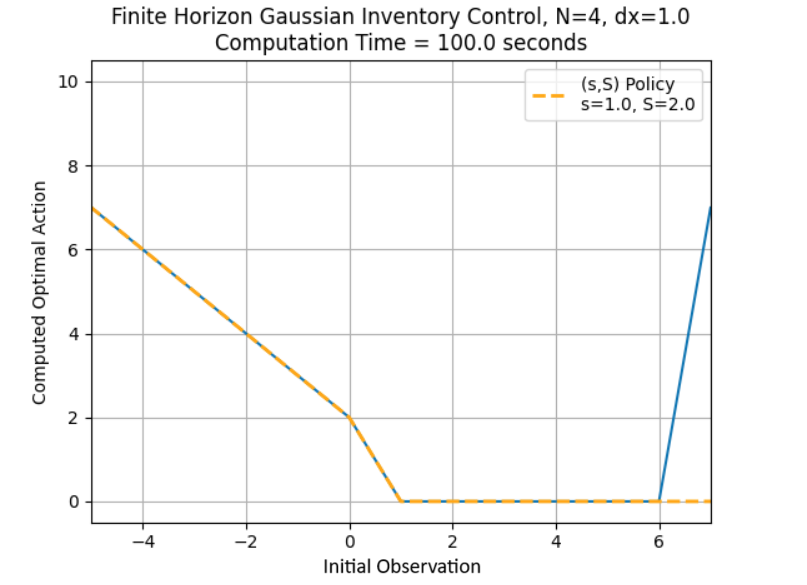} &   \includegraphics[width=65mm]{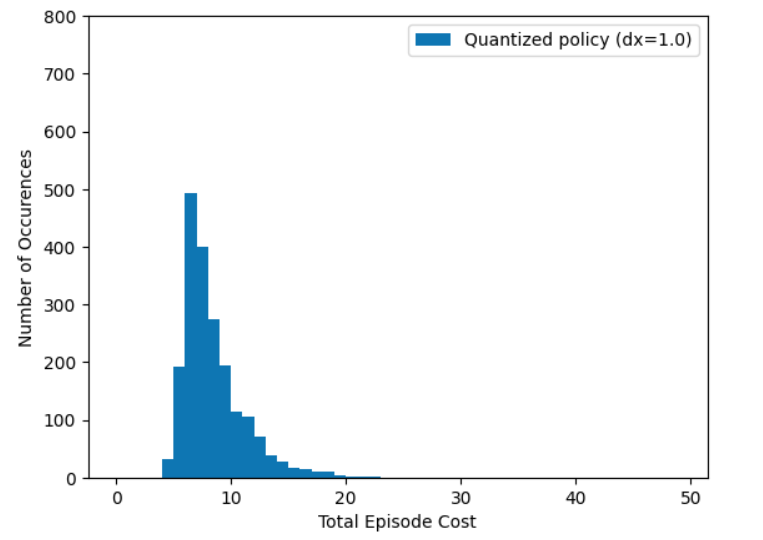} \\
 \includegraphics[width=65mm]{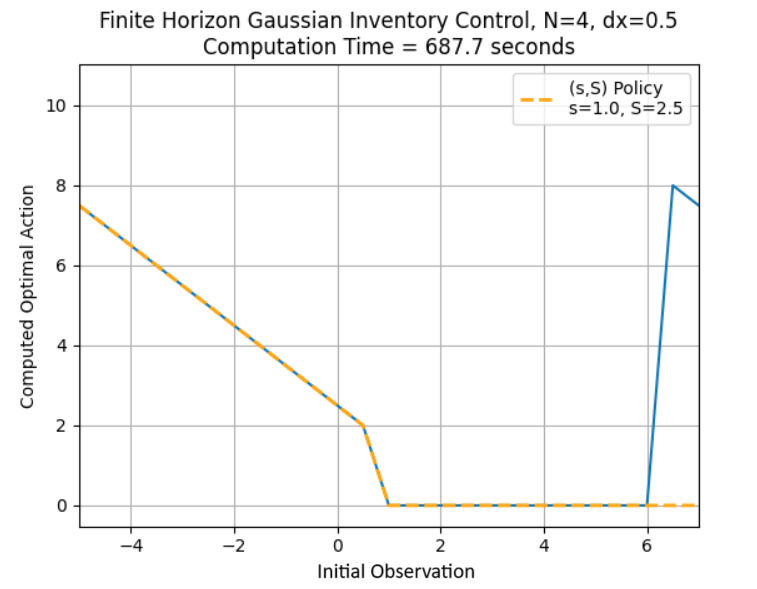} &   \includegraphics[width=65mm]{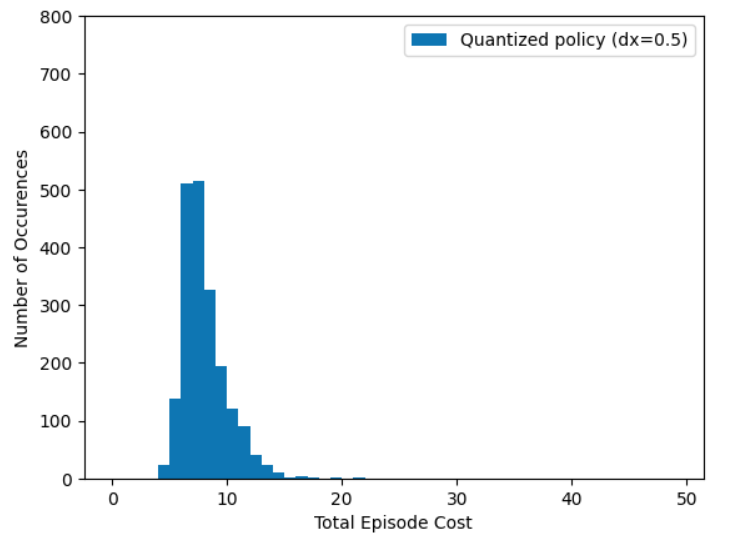} \\{\includegraphics[width=65mm]{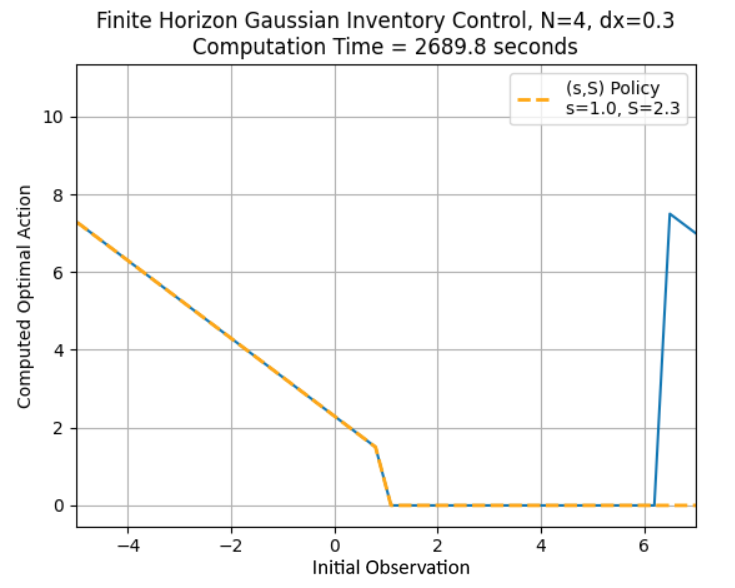} } & {\includegraphics[width=65mm]{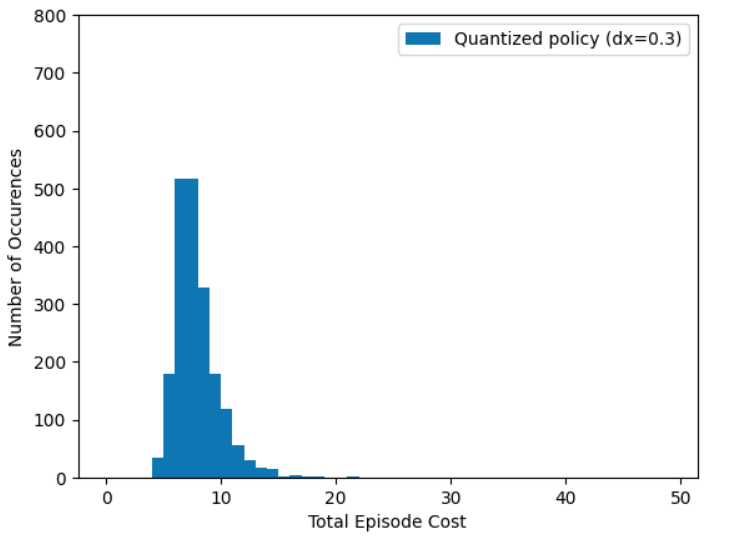} } \\
\end{tabular}
\caption{Structure and performance of policies derived from quantization (Gaussian case).}
    \label{fig:quantized}
\end{figure}

\begin{figure}

    \centering
    \begin{tabular}{cc}

        {\includegraphics[width=65mm]{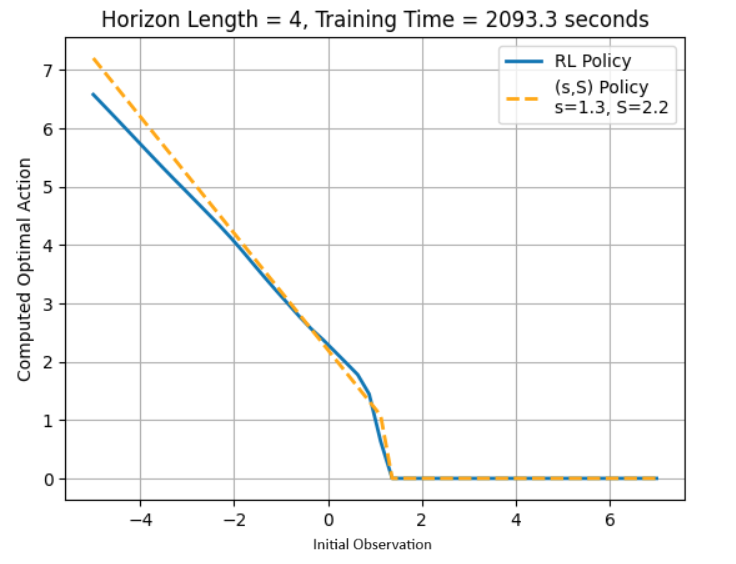} } & {\includegraphics[width=65mm]{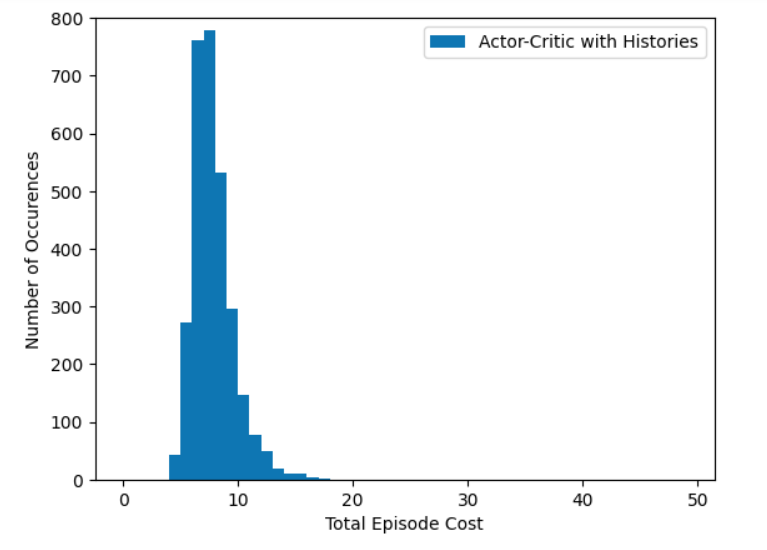} }
     \end{tabular}
     \caption{Structure and performance of policy computed via DDPG applied to histories (Gaussian case).}
     \label{fig:ddpg_histories}
\end{figure}

\begin{figure}

    \centering
    \begin{tabular}{cc}
         {\includegraphics[width=65mm]{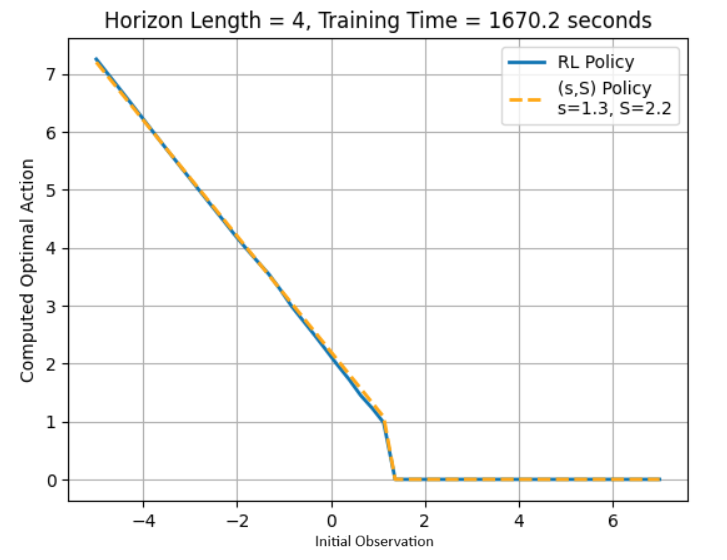} } & {\includegraphics[width=65mm]{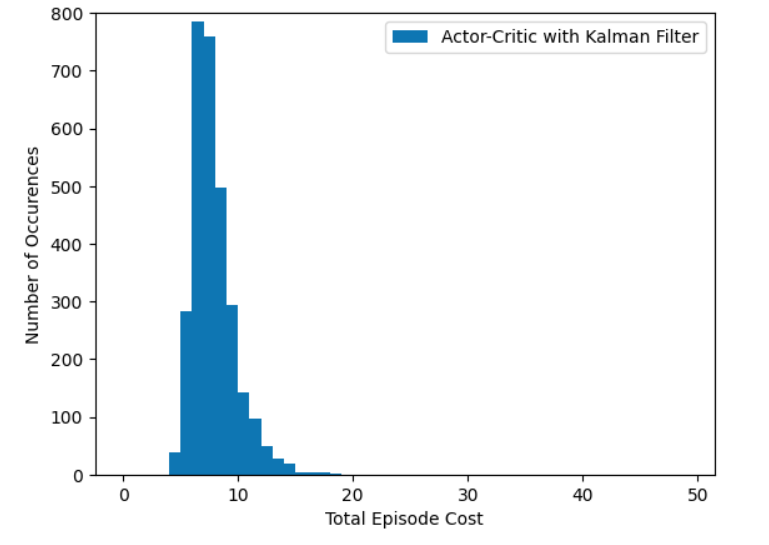} }
     \end{tabular}
     \caption{Structure and performance of policy computed via DDPG applied to beliefs (Gaussian case).}
         \label{fig:ddpg_beliefs}
\end{figure}

\begin{figure}
    \centering
    \begin{tabular}{cc}
         {\includegraphics[width=65mm]{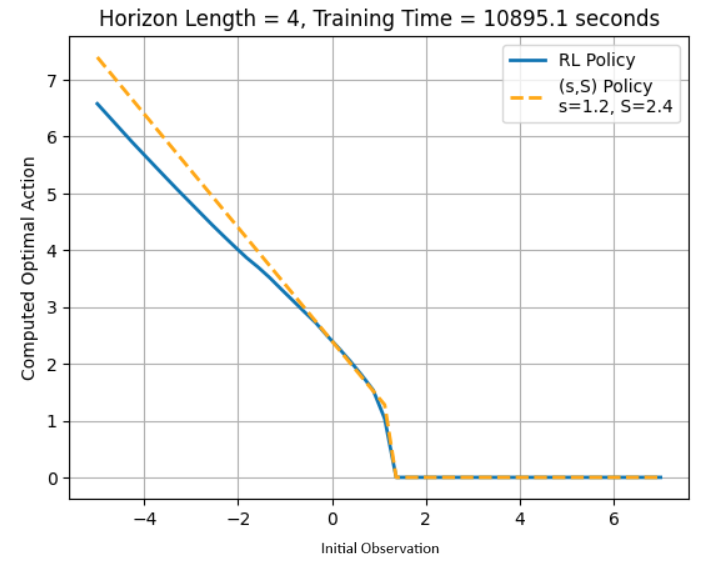} } & {\includegraphics[width=65mm]{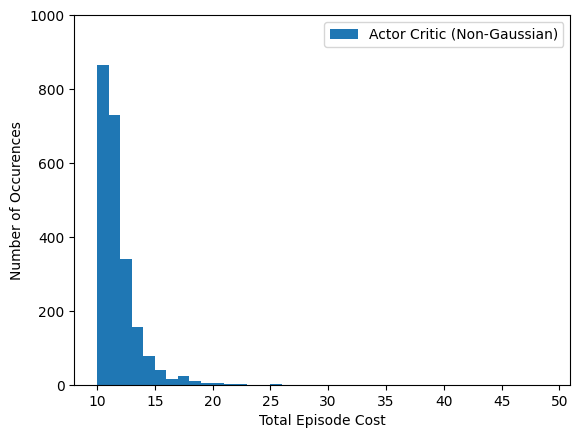} }
     \end{tabular}
    \caption{Structure and performance of policy computed via DDPG applied to histories (Non-Gaussian case).}
    \label{fig:ddpg_non_gaussian}
\end{figure}

The structure of the computed policies, and the distribution of their performance across 3000 episodes, are presented in Figures \ref{fig:quantized}, \ref{fig:ddpg_histories}, \ref{fig:ddpg_beliefs}, and \ref{fig:ddpg_non_gaussian}. To illustrate the structure of the computed policies, we plotted the selected action as a function of the observation obtained at time $t = 0$. The plots in the left-hand column of Figure~\ref{fig:quantized} show that at time $t = 0$, the computed policies behave similarly to an $(s, S)$ policy. A similar phenomenon is observed for the policies computed using DDPG in both the Gaussian (Figures \ref{fig:ddpg_histories} and \ref{fig:ddpg_beliefs}) and non-Gaussian (Figure \ref{fig:ddpg_non_gaussian}). In fact, in the Gaussian case we prove (\ref{secA3}) that there exists a time-dependent $(s, S)$ policy that is optimal.


\section{Conclusions}
This paper studies the computation of optimal policies for a partially observable single-commodity periodic review inventory control problem with backorders and setup costs.  The goal is to maximize the expected total discounted cost incurred over a finite planning horizon. This problem can be reduced to a belief MDP where the state space is the space of all probability distributions on the real line, so that every state is an infinite-dimensional object. To avoid this complication, we applied the DDPG algorithm to an MDP where the states are the observed histories of the original problem. This approach becomes more difficult to implement as the length of the planning horizon increases, because of the growth in dimensionality of the space of histories.

We also investigated the Gaussian version of this problem with normally distributed initial states, demands, and observation noise. We showed that the belief MDP for this problem can be reduced to an MDP corresponding to a single-commodity periodic review inventory control problem with time-dependent demands and holding costs.  The state space for this MDP is the real line, where each state is the mean of the posterior distribution of the true inventory level. This formulation does not suffer from increases in the time horizon.  We also applied DDPG to this problem, which produced policies that behave and perform similarly to estimates of optimal policies obtained via discretizations. This is consistent with our theoretical analysis for the Gaussian case, which indicates that there are optimal time-dependent $(s, S)$ policies for the original problem.

\appendix

\renewcommand{\appendixname}{Appendix}
\renewcommand{\thesection}{\appendixname~\Alph{section}}

\section {POMDPs for Stochastic Filtering Problems} \label{secA1}
A periodic review inventory control problem with noisy observations  is a particular case of the following nonlinear filtering problem with $\Xc:=\R$ and  $\Hc:=\R$ being the spaces of state and observation noises, for $t = 0,1, \dots,$ 
\begin{subequations}\label{eq:filters:model}
    \begin{align}
        x_{t+1} &= F(x_t, a_t, \xi_t), && x_t \in \Xb, \quad a_t \in \Ab, \quad \xi_t \in \Xc, \label{eq:filters:transitions} \\
        y_{t+1} &= G(a_t, x_{t+1}, \eta_{t+1}), && {\color{black}a_t \in \Ab, \quad x_{t+1} \in \Xb,} 
        \quad \eta_{t+1} \in \Hc, \label{eq:filters:observations}
    \end{align}
\end{subequations}
where $F : \Xb \times \Ab \times \Xc \to \Xb$ is the Borel-measurable functions representing the transition dynamics of the system, and $G : \Ab \times \Xb \times \Hc \to \Yb$ is the Borel-measurable function representing the observations. In addition, there is the probability distribution $p_0 \in \Pb(\Xb)$ of the initial state $x_0,$ and $y_0= G_0(x_0,\eta_0),$  where $G_0 : \Xb \times \Hc \to \Yb$ a Borel-measurable function defining the initial observation, and
$\set{\xi_t}_{t=0}^{\infty}$ and $\set{\eta_t}_{t=0}^{\infty}$ are sequences of independent and identically distributed (iid) random variables with distributions $\mu \in \Pb(\Xc)$ and $\nu \in \Pb(\Hc),$ respectively. It is also assumed that the sequence $\set{x_0, \xi_0, \eta_0, \xi_1, \eta_1, \dots}$ is mutually independent. The process evolves as follows. {\color{black} The initial hidden state $x_0$ has the distribution $p_0,$ and the initial observation is $y_0 = G(x_0, \eta_0).$  At each epoch $t = 0, 1, \dots,$  the decision maker, who  knows the initial state distribution $p_0,$ the previous and current observations, and the previously selected actions $p_0,y_0,a_0,\ldots,y_{t-1}a_{t-1}y_t,$ selects the action $a_t \in \Ab.$} The next state $x_{t+1}$ and observation $y_{t+1}$ are determined by equations \eqref{eq:filters:model}. 
In particular, in the described above the inventory control problem, for $t=0,1,\ldots,$  $\xi_t=D_t,$ \[F(x_t, a_t,\xi_t)=L( x_t + a_t - \xi_t),\qquad G(a_t, x_{t+1}, \eta_{t+1})= x_{t+1}+ \eta_{t+1},\] and $G_0(x_0, \eta_0)= x_0+ \eta_0.$

The evolution of a system defined by stochastic equations such as~\eqref{eq:filters:model} can also be represented in terms of stochastic kernels, and the corresponding model is a POMDP.  For POMDPs, transitions of states are defined by a stochastic kernel $\Tc$ on $\Xb$ given $\Xb\times\Ab$ called the transition kernel or transition probabilities, and observations are defined by a stochastic kernel $Q$ on $\Yb$ given $\Ab\times\Xb$ called the observation kernel or observation probabilities.  For POMDPs, the state $x_{t+1}$ is defined by the distribution $\Tc(\:\cdot\:|x_t,a_t)$ instead of by equality~\eqref{eq:filters:transitions}, and the observation $y_{t+1}$ is defined by the distribution $Q(\:\cdot\:|a_t,x_{t+1})$  instead of by equality~\eqref{eq:filters:observations}. In addition, the initial distribution of observation $y_0$ is defined by a stochastic kernel $Q_0$ on $Y$ given $X,$ that is, $y_0$ is defined by the distribution $Q(\:\cdot\:|x_0);$ see \cite[Chapter 4]{HL} or \cite{feinberg_partially_2016} for details.

The model definitions based on stochastic equations and on POMDPs are equivalent.  Indeed, for the functions $F$ and $G$ from~\eqref{eq:filters:model},
\begin{subequations} \label{eq:filters:model_kernels}
    \begin{align}
        \Tc(B|x, a) &= \int_{\Xc} \mathbf{1}\{F(x, a, \xi) \in B\} \ \mu(d\xi), && B \in \Bc (\Xb), \quad x \in \Xb, \quad a \in \Ab, \label{eq:filters:transition_kernel} \\
        Q(C|a, x) & =\int_{\Hc}\mathbf{1}\{G(a, x, \eta) \in C\} \ \nu(d\eta), && C \in \Bc (\Yb), \quad a \in \Ab, \quad x \in \Xb, \label{eq:filters:observation_kernel}
    \end{align}
\end{subequations}
are {\color{black}stochastic kernels; see~\citet[Proposition C.2]{HL}} or 
\citet[p. 190]{bertsekas_stochastic_1996}. Of course, the function $G_0$ also defines the stochastic kernel $ Q_0(C|a, x)  =\int_{\Hc}\mathbf{1}\{G_0( x, \eta) \in C\}\nu(d\eta),$ where $ C \in \Bc (\Yb),$ $a \in \Ab,$ and $x \in \Xb.$
The reduction of a POMDP to the stochastic system defined by equations~\eqref{eq:filters:model}, with $\Xc=[0,1],$ $\Hc=[0,1],$ and $\mu$ and $\eta$ are uniform distributions,  follows from \citet[Lemma F]{aumann_1964_mixed}; see \cite{FIKK} for details. However, we do not assume here that $\Xc=[0,1],$ $\Hc=[0,1],$ and $\mu$ and $\eta$ are  {\color{black}uniform distributions} because sometimes it can be more convenient to consider other distributions including normal distributions considered in this paper.

 \section {Belief-MDPs for Stochastic Filtering Problems}  \label{secA2}
 The classic approach to solving a POMDP is based on its reduction to a completely observable MDP, whose states are probability distributions of the states of the system, sometimes called belief distributions, and this MDP is sometimes called a filter.  Weak continuity of the transition probabilities of the filter, also called weak continuity of the filter, is an important property.  If one-step costs are $\Kb$-inf-compact and bounded below, then weak continuity of the filter implies the existence of optimal policies for problems with expected total costs, the validity of optimality equations, and convergence of value iterations; see \cite{feinberg_partially_2016, feinberg_markov_2022, FIKK} for details.


Let us consider the standard reduction of a POMDP with a transition kernel $\Tc$ and observation kernel $Q$ to a filtered MDP with complete observations and with the transition kernel $q$ to be defined in~\eqref{eq:filters:filter_kernel}. For more details on this reduction, see~\citet[Section~10.3]{bertsekas_stochastic_1996},~\citet[Chapter~8]{dynkin_1979_controlled},~\citet[Chapter~4]{HL},~\citet{rhenius_1974}, or~\citet{yushkevich_1976_reduction}. As usual, we introduce the stochastic kernel on $\Xb \times \Yb$ given $\Xb \times \Ab$ defined by the equation
\begin{align}\label{eq:filters:joint_kernel}
    P(B \times C|x,a) &\defeq
    \int_{B} Q(C|a,x') \ \Tc(dx'|x,a), &&
    B \in \Bc(\Xb), \quad
    C \in \Bc(\Yb), \quad
    x \in \Xb, \quad
    a \in \Ab,
\end{align}
which {\color{black}represents the joint probability of the next state-observation pair} 
$(x', y')$ given the state-action pair $(x, a).$  
Consider a measure $z \in \Pb(\Xb)$ representing a belief distribution for the current state $x.$ The joint distribution of the next state-observation pair $(x', y')$ conditioned by the measure $z$ and action $a \in \Ab$ is given as the stochastic kernel
\begin{align} \label{eq:filters:joint_belief_kernel}
    R(B \times C | z,a) &\defeq
    \int_{\Xb} P(B \times C | x,a) \ z(dx), &&
    B \in \Bc(\Xb), \quad
    C \in \Bc(\Yb), \quad
    z \in \Pb(\Xb), \quad
    a \in \Ab,
\end{align}
with marginal distribution for $y'$ given by $R'(C|z,a) \defeq R(\Xb \times C | z,a).$ According to~\citet[Proposition 7.27]{bertsekas_stochastic_1996} there exists a stochastic kernel $H$ {\color{black}on $\Xb$ given $\Pb(\Xb) \times \Ab\times\Yb$} 
such that
\begin{align} \label{eq:filters:joint_belief_conditional_representation}
    R(B \times C | z, a) &=
    \int_C H(B|z,a,y) \ R'(dy|z,a), &&
    B \in \Bc(\Xb), \quad
    C \in \Bc(\Yb), \quad
    z \in \Pb(\Xb), \quad
    a \in \Ab.
\end{align}
Alternatively, the stochastic kernel $H$ defines a Borel measurable mapping $(z,a,y') \mapsto H(z,a,y') \in \Pb(\Xb)$ that describes the evolution of the distributions $z$ by $z'= H(z, a, y').$ We remark that the mapping $H$ is a generalization of filters commonly found in applications. In particular, for linear state space models with {\color{black}a Gaussian initial state distribution $p_0$ and with Gaussian noises}, the belief $z$ is a Gaussian distribution parameterized by {\color{black}a mean and a covariance.} 
The Kalman filter, which computes the next distribution $z'$ from $z,$ $a$ and the next observation $y',$ computes $H$ explicitly in this case.  Transition dynamics from the current state distribution $z$ to the next state distribution $z'$ is defined by the stochastic kernel $q$ defined as
\begin{align} \label{eq:filters:filter_kernel}
    q(D |z, a) &\defeq
   \int_{\Yb} \mathbf{1}\{H(z, a, y) \in D\} \ R^{\prime}(dy |z, a), &&
    D \in \Bc(\Pb(\Xb)), \quad
    z \in \Pb(\Xb), \quad
    a \in {\color{black}\Ab;}
\end{align}
see, e.g., \cite{bertsekas_stochastic_1996, dynkin_1979_controlled, feinberg_partially_2016, HL, rhenius_1974, yushkevich_1976_reduction}.
The stochastic kernel $q$ is the transition probability for the filter.

Let us consider again the model defined by equations \eqref{eq:filters:model}.  In view of~\eqref{eq:filters:model_kernels},
\begin{equation} \label{eq:filters:joint_kernel_representation}
    P(B \times C | x,a) =\int_{\Xc}\int_{\Hc} \mathbf{1}\{G(a,F(x,a,\xi),\eta) \in C\} \mathbf{1}\{F(x,a,\xi)\in B\} \ \nu(d\eta) \ \mu(d\xi).
\end{equation}
We observe that formula~\eqref{eq:filters:joint_kernel_representation} employs the assumption of joint independence of the noise variables $\xi$ and $\nu$ mentioned above.

For the inventory control problem,
\begin{equation} \label{eq:filters:joint_kernel_inv}
    P(B \times C | x,a) =\int_{-\infty}^{+\infty}\int_{-\infty}^{+\infty} \mathbf{1}\{L( x + a - d)+w \in C\} \mathbf{1}\{L( x + a - d)\in B\} \ dF_\eta(w) \ dF_D(d).
\end{equation}
The Belief MDP is the MDP with the state space $\Pb(\X),$ action sets $\Ab,$ transition probability $q,$ and one-step costs $C,$ where $C(z,a)=\int_{-\infty}^{+\infty} c(x,a)z(dx)$  for $z\in\Pb(\X)$ ans $a\in\Ab.$  In our case, $\X=\R$ and  $\Pb(\X)$   is the metric space of probability measures on a real line  endowed with the L\'evy-Prokhorov metric or any other metric defining weak convergence.

Let us consider expected discounted costs.  Let $V_{\alpha,t}(z)$ be the $t$-horizon value function for the belief MDP, where $t=1,2,\ldots,$ and $V_\alpha(z)$ is the value function for the infinite horizon. According to \cite[Corollary 7.6]{FIKK}, the finite-horizon value functions satisfy the optimality equation
\[ V_{\alpha,t+1}(z)=\min_{a\ge 0} \{C(z,a)+\int_{\Pb(\X)} V_{\alpha,t}(z')q(dz'|z,a)\}, \qquad t=0,1,\ldots,\ z\in\Pb(\X),\]
where $ V_{\alpha,0}(z)\equiv 0.$ In addition, $V_{\alpha,t}\uparrow V_\alpha$ as $t\to\infty,$ and the infinite-horizon value function also satisfies the optimality equation
\[ V_\alpha (z)=\min_{a\ge 0} \{C(z,a)+\int_{\Pb(\X)} V_\alpha (z')q(dz'|z,a)\},\qquad z\in\Pb(\X).\]
The finite-horizon optimality equation defines Markov optimal policies finite-horizon belief MDPs, while the infinite-horizon optimality equation defines  deterministic optimal policies for infinite-horizon belief MDPs.  These optimal policies define optimal policies for the original problem by selecting at step $t$ the optimal action for the belief MDP with the state equal to posterior distribution of states after $t$th observation becomes known.

These result provide important insights, but they lead to inefficient computations and also require to deal with infinite-dimensional spaces if the original problem has an infinite state space.  One of the ways to avoid dealing with spaces of probability measures is to consider the vectors of previous observations $h_t=(y_0,a_0,\ldots, y_{t-1},a_{t-1}y_t).$  If an action $a_t$ is selected, then  the next step will become $(h_t,a_t,y_{t+1})$ with the probability $R'(\cdot|z_t,a_t),$ where $z_t$ is the belief distribution defined by the initial state distribution $p_0$ and by the history $h_t.$ Of course, computing these probabilities is a nontrivial issue, and therefore we use the RL methods that do not require the knowledge of transition probabilities. The important factor is that these probabilities exist.

\section {Gaussian Inventory Control Model with Backorders} \label{secA3}

For a normal random variable $X,$ we write $X\sim\mathcal{N}(\bar{X},\sigma_{X}^2),$ where $\bar{X}:= \mathsf{E}[X]$ and $\sigma_{X}^2={\rm Var}(X),$ where $\sigma_{X}\ge0$ is the standard deviation of $X.$ If $(X,Y)$ is a bivariate normal random vector, then
\begin{equation}\label{eqXfromY} X|Y=y \sim \mathcal{N}\left(\bar{X}+\frac{{\rm Cov}(X,Y)}{\sigma_Y^2}(y-\bar{Y}), \sigma_X^2-\frac{{\rm Cov}^2(X,Y)}{\sigma_Y^2}\right).\end{equation}

Let us consider the Gaussian version of the inventory control problem with back orders.  This means that  $x_0\sim \mathcal{N}(\bar{x}_0,\sigma_{x_o}^2),$ $D_t\sim \mathcal{N}(\bar{D}_t,\sigma_{D_t}^2),$ and $\eta_t\sim \mathcal{N}(0,\sigma_{\eta_t}^2),$ $t=0,1,\ldots,$ are independent normal random variables.  We also have that
\begin{equation} \label{eqinvmicolll} x_{t+1}=x_t+a_t-D_t\quad {\rm and}\quad y_t=x_t+\eta_t, \qquad t=0,1,\ldots. \end{equation}
This is a linear system of equations for normal random variable, and $(x_t,y_t)$ are bivariate normal random vectors.  Let $z_t$ be posterior distributions of state $x_t$ at time $t=0,1,\ldots$ after observations $h_t:=y_0,a_0,y_1,\ldots,y_{t-1},a_{t-1},y_t $ are known.  Let $X_t$ be random variables with distributions $z_t.$ Since $\mathsf{E}[\eta_0]=0,$  we have that $\bar{y}_0=\bar{x}_0,$ $\sigma_{y_0}^2=\sigma_{x_0}^2+\sigma_{\eta}^2,$  and ${\rm Cov}(x_0,y_0)=\sigma_{x_0}^2.$ Formula \eqref{eqXfromY} applied to $X=x_0$ and $Y=y_0$ implies \begin{equation}\label{eqnormdistforxp}X_0:=X|Y=y_0\sim\mathcal{N}\left(\frac{\sigma_\eta^2}{\sigma_{x_0}^2+\sigma_\eta^2}\bar{x}_0+\frac{\sigma_{x_0}^2}{\sigma_{x_0}^2+\sigma_\eta^2}y_0,\frac{\sigma_{x_0}^2\sigma_\eta^2}{\sigma_{x_0}^2+\sigma_\eta^2}\right)=z_0\end{equation} is a normal random variable. Formula \eqref{eqnormdistforxp} implies that $\sigma_0^2:={\rm Var}(X_0)=\frac{\sigma_{x_0}^2\sigma_\eta^2}{\sigma_{x_0}^2+\sigma_\eta^2},$ and this variance does not depend on $y_0.$  Since $y_0$ is random,  \[\bar{X}_0:=\frac{\sigma_\eta^2}{\sigma_{x_0}^2+\sigma_\eta^2}\bar{x}_0+\frac{\sigma_{x_0}^2}{\sigma_{x_0}^2+\sigma_\eta^2}y_0\]
is the normal random variable with \begin{equation}\label{eqharfor0}\mathsf{E}[\bar{X}_0]=\bar{x}_0\quad{\rm and}\quad {\rm Var}(\bar{X}_0)={\rm Var}\left(\frac{\sigma_{x_0}^2}{\sigma_{x_0}^2+\sigma_\eta^2}y_0\right)=\left(\frac{\sigma_{x_0}^2}{\sigma_{x_0}^2+\sigma_\eta^2}\right)^2 \sigma_{y_0}^2=\frac{\sigma_{x_0}^4}{\sigma_{x_0}^2+\sigma_\eta^2}.\end{equation}

Let  $z_t=\mathcal{N}(\bar{X}_t,\sigma_t^2)$ be the posterior distribution of the inventory level at time  $t=0,1,\ldots$ given a history $h_t,$ and let $X_t\sim \mathcal{N}(\bar{X_t},\sigma_t^2)$ be the posterior estimation of the state at time $t.$ In particular, $X_0$ is defined in  \eqref{eqnormdistforxp}. Then  at time $(t+1),$ given the history $h_t,a_t,$ the inventory level is
 \[X_t+a_t-D_t\sim\mathcal{N}(\bar{X}_t+a_t-\bar{D},\sigma_t^2+\sigma_D^2),\]   the observation is \[Y_{t+1}=X_t+a_t-D_t+\eta_{t+1}\sim \mathcal{N}(\bar{X}_t-\bar{D},\sigma_t^2+\sigma_D^2+\sigma_\eta^2),\] and therefore $(X_t+a_t-D_t,Y_{t+1}) $ is a bivariate normal random variable with
\[
{\rm Cov}(X_t+a_t-D_t,Y_{t+1})=\sigma_t^2+\sigma_D^2.\]
At time (t+1), the prior estimation of the state is $X_t+a_t-D_t,$ and the posterior estimation of the state is $X_{t+1}.$  Formula \eqref{eqXfromY} applied to $X=X_t+a_t-D_t$ and $Y=X+\eta_{t+1}$ implies
\begin{equation}\label{eqxtpp}
X_{t+1}=X|Y\sim
\mathcal{N}\left(\frac{\sigma_\eta^2}{\sigma_t^2+\sigma_D^2+\sigma_\eta^2}(\bar{X}_t+a _t-\bar{D})+\frac{\sigma_t^2+\sigma_D^2}{\sigma_t^2+\sigma_D^2+\sigma_\eta^2}Y,\frac{(\sigma_t^2+\sigma_D^2)\sigma_\eta^2}{\sigma_t^2+\sigma_D^2+\sigma_\eta^2}\right)=z_{t+1},\end{equation} and, as defined above, $Y~\sim\mathcal{N}(\bar{X}_t+a_t-\bar{D}_t\sigma_t^+\sigma_D^2+\sigma_\eta^2)$ is the prior distribution of the observation   $y_{t+1}.$ Thus
 \begin{equation}
 \bar{X}_{t+1}:=\frac{\sigma_\eta^2}{\sigma_t^2+\sigma_D^2+\sigma_\eta^2}(\bar{X}_t+a_t-\bar{D})+\frac{\sigma_t^2+\sigma_D^2}{\sigma_t^2+\sigma_D^2+\sigma_\eta^2}Y\end{equation}
 is the normal random variable with
\begin{equation}\label{echarfortplus1}\mathsf{E}[\bar{X}_{t+1}]=\bar{X}_t+a_t-\bar{D}\quad{\rm and}\quad {\rm Var}(\bar{X}_{t+1})={\rm Var}(\mathsf{E}[X_t+a_t-D_t|y_{t+1})=\frac{(\sigma_t^2+\sigma_D^2)^2}{\sigma_t^2+\sigma_D^2+\sigma_\eta^2},\end{equation}
where \eqref{echarfortplus1} is similar to \eqref{eqharfor0}.

Thus,
\begin{equation}\label{eqnormdfxpt}
 \bar{X}_{t+1} ~\sim \mathcal{N}(\bar{X}_t+a_t-\bar{D},\frac{(\sigma_t^2+\sigma_D^2)^2}{\sigma_t^2+\sigma_D^2+\sigma_\eta^2}).
\end{equation}
Formula \eqref{eqnormdfxpt} implies that $\bar{X}_{t+1}$ depends on  $y_0,a_0,\ldots,y_{t-1},a_{t-1},y_t$ only via the values of $\bar{X}_t,a_t,$ and $y_t.$  Therefore,
\begin{equation}\label{eqnormdfxptt}
 \bar{X}_{t+1}|\bar{X}_t ~\sim \mathcal{N}(\bar{X}_t+a_t-\bar{D},\frac{(\sigma_t^2+\sigma_D^2)^2}{\sigma_t^2+\sigma_D^2+\sigma_\eta^2}),
\end{equation}
and this formula can be rewritten as
\begin{equation}\label{eqtranfdsmdp}
 \bar{X}_{t+1}=\bar{X}_t+a_t-D^*_t,\qquad t=0,1,\ldots,
\end{equation}
where \begin{equation}\label{eqdefDSTAR}D_t^*\sim N(\bar{D},\frac{(\sigma_t^2+\sigma_D^2)^2}{\sigma_t^2+\sigma_D^2+\sigma_\eta^2})\end{equation}
are independent random variables.
Formula \eqref{eqtranfdsmdp} indicates that expectations for beliefs form a periodic review inventory control problem with back orders and with variable demand distributions.

We recall that $\sigma_0^2=\frac{\sigma_{x_0}^2\sigma_\eta^2}{\sigma_{x_0}^2+\sigma_\eta^2}.$  If $\sigma_{x_0}=\sigma_\eta=0,$ we deal with the the trivial case when there is no demand and there is no noise. Let us exclude this case.  Then the value of $\sigma_0$ is well-defined. Formula \eqref{eqxtpp} implies
\begin{equation}\label{eqstp1ef}
\sigma_{t+1}^2=\frac{(\sigma_t^2+\sigma_D^2)\sigma_\eta^2}{\sigma_t^2+\sigma_D^2+\sigma_\eta^2},\qquad t=0,1,\ldots.
\end{equation}
The sequence $\{\sigma_t^2\}_{t=0}^\infty$ has the limit
\[ \sigma_*^2:=\lim_{t\to\infty}\sigma_t^2=\frac{1}{2}\sigma_D^2(\sqrt{1+\frac{4\sigma_\eta^2}{\sigma_D^2}}-1)\]
if $\sigma_D\ne 0.$ This follows from the Banach fixed point theorem applied to the function $x\mapsto \frac{(x+\sigma_D^2)\sigma_\eta^2}{x+\sigma_D^2+\sigma_\eta^2}.$ If  $\sigma_D= 0$ or  $\sigma_{\eta}= 0,$ then $\sigma_*=0.$  If $\sigma_\eta^2=2\sigma_D^2,$ then $\sigma_*=\sigma_D.$

    Since all beliefs are Gaussian, the belief MDP can be simplified by parameterizing normal distributions.  So, each Gaussian belief $z$ can be presented by a pair $(\bar{x},\sigma^2),$
which are the mean and the variance of the normal distribution corresponding to a Gaussian probability $z$ on $\R.$  In other words, for Gaussian distributions, the state space of all probability measures on $\R$ can be represented with the two-dimensional space $\R\times\R_{\ge 0},$ where $\bar{x}\in\R$ is the mean and $\sigma^2\in \R_{\ge 0}$ is the variance.

In Appendix A, we denoted by $x_t$ the states after observations $y_t$ were known, and $x_t\sim z_t,$ where $z_t$ is the belief distribution at time $t.$
This $x_t$ correspond to the notation $X_t$ in this section. The notation $\bar{x}_t$ we use below corresponds to $\bar{x}_t.$  

So, the belief MDP can be replaced with its parametrized version with the state space $\R\times\R_{\ge 0}$ and  action sets $\A=\R_{\ge 0}.$  If at time $t=0,1,\ldots$ an action $a_t\in\A$ is chosen at a state $(\bar{x}_t,\sigma_t^2)\in \R\times\R_{\ge 0},$ then the next state is $(\bar{x}_{t+1},\sigma_{t+1}^2), $ where $\bar{x}_{t+1}$ has the normal distribution defined in \eqref{eqnormdfxpt} or\eqref{eqnormdfxptt}, and $\sigma_{t+1}^2$ is defined in \eqref{eqstp1ef}.
In the functional form, we can write the changes of $\bar{x}_t$ as
\begin{equation}\label{eqsystdy}
    \bar{x}_{t+1}=\bar{x}_t+a_t-D_t^*,\qquad t=0,1,\ldots,
\end{equation}
where $D_t^*$ are defined in \eqref{eqdefDSTAR}, and \eqref{eqsystdy} is an MDP version of \eqref{eqtranfdsmdp}.
The one-step cost is $c^*((\bar{x}_t,\sigma_t^2),a_t)=\mathsf{E}[c(\xi_t,a_t)]:=\hat{c}(a_t)+\mathsf{E}[h(\xi_t+a_t-D_t)],$ where $\xi_t\sim\mathcal{N}(\bar{x}_t,\sigma_t^2),$ and $\xi_t$ and $D_t$ are independent.  In view of \eqref{eqharfor0}, it is  natural to consider initial states in the form of $(\bar{x}_0,\frac{\sigma_{x_0}^4}{\sigma_{x_0}^2+\sigma_{\eta}^2}).$

This problem satisfies \cite[Assumption (W*)]{FKZ} for MDPs with weakly continuous transition probabilities.  Therefore, there exists Markov optimal policies for problems with finite horizons, deterministic optimal policies for infinite-horizon problems.  In addition, optimality equations hold for this model, they satisfy optimality equations, and these equations identify optimal policies.  Value iteration algorithms converge for this problem.

The special property of this  problem is that the $(\sigma_t)_{t=0}^\infty$ is the deterministic sequence which does not depend on observations and selected actions.  Therefore, this problem  with two-dimensional states $(\bar{x},\sigma^2)$ can be reduced to the problem with one-dimensional states $\bar{x}$ and with transition probabilities and costs depending on the time parameter. In other words, the DM  knows only expectations of beliefs and time parameters.

Let us consider an MDP with the state space $\R$ and action sets $\A=\R_{\ge 0}.$   If at a state $\bar{x}_t$ an action $a_t$ is chosen at time $t=0,1,\ldots,$
then the system moves to the next $\bar{x}_{t+1}$ according to equation  \eqref{eqsystdy}.  The one-step cost is  $c_t^*(\bar{x}_t,a_t)=c^*((\bar{x}_t,\sigma_t^2),a_t).$  We would like to rewrite $c_t^*(\bar{x}_t,a_t)$ in the form
\begin{equation}\label{eqnewholc}
c^*_t(\bar{x}_t,a_t)=\hat{c}(a_t)+\mathsf{E}h^*_t(\bar{x}_{t}+a_t-D^*_t),
\end{equation}
where $h_t^*:\R\mapsto\R\times \R_{>0}$ is convex and $h_t^*(x)\to\infty$ as $|x|\to\infty$ for $t=0,1,\ldots.$  In this MDP, transition probabilities and costs are time-dependent.

We observe that $\xi_t+a_t-D_t\sim \mathcal{N}(\bar{x}_t+a_t-\bar{D},\sigma_t^2+\sigma_D^2).$ while $x_t+a_t-D^*_t\sim \mathcal{N}(\bar{x}_t+a_t-\bar{D},\sim \mathcal{N}(\bar{D},\frac{(\sigma_t^2+\sigma_D^2)^2}{\sigma_t^2+\sigma_D^2+\sigma_\eta^2}).$ We also observe that
\[
\sigma_{t+1}^2:=\sigma_t^2+\sigma_D^2-\frac{(\sigma_t^2+\sigma_D^2)^2}{\sigma_t^2+\sigma_D^2+\sigma_\eta^2}=\frac{\sigma_\eta^2(\sigma_t^2+\sigma_D^2)}{\sigma_t^2+\sigma_D^2+\sigma_\eta^2}.
\]
Let us consider a sequence of independent random variables $\{Z_t\}_{t=0}^\infty$ with $Z_t\sim \mathcal{N}(0,\sigma_{t+1}^2)$ and define
\[
h^*_t(x):=\mathsf{E}[h(x+Z_t)].
\]
This function is nonnegative, concave, and, in view of the monotone convergence theorem, $h_t^*(x)\to\infty$ as $|x|\to\infty.$ Since $h(\xi_t+a_t-D_t)\sim h_t^*(x_t+a_t-D_t^*),$  formula \eqref{eqnewholc} holds.  Thus, the defined MDP with time-dependent transition probabilities and costs describes an inventory control problem with time-dependent demand and holding costs.  The states in this MDP are expectations of prior distributions of inventory levels in the original inventory control problem.   This explains why the optimal policies for expected prior estimations have the $(s_t,S_t)$ structure, which was originally picked up by the RL algorithm.

Here we would like to recall  the structure of optimal policies for periodic review problems with complete observations and holding costs and demands depending on the time parameter $t.$ That is, the holding cost at time $ (t+1)$ is $h_t(x_{t+1})$ and $\{D^*_t\}_{t=0}^\infty$ are independent nonnegative random demands, $D^*_t\sim F_{D^*_t}.$ We assume that the functions $h_t$ are nonnegative and convex,  $h_t(x)\to \infty$ as $t\to\infty$ and $\mathsf{E}D_t(y-D^*_t)<\infty$ for all $y\in\R$ and $t=0,1,\ldots.$ Analysis of such problems and relevant models and results can be found in the literature; see, e.g., \cite{heymansobel, simchilevi, bensinv, BCST, feinberglewis, feinbergliang}.

A Markov policy $\phi$ is called an $(s_t,S_t)$-policy if for each $t=0,1,\ldots$ there exist numbers $s_t$ and $S_t,$ where $s_t <S_t,$ such that at time $t$ the policy $\phi$ does not order inventory, if $x_t\ge 0,$ and orders up to the level $S_t$ if $x_t<s_t.$ The following condition introduced by  \cite{VeinottW} is sufficient for the existence of optimal $(s_t,S_t)$-policies for problems with expected total discounted costs. For each $t=0,1,\ldots$ there exist $z_t,y_t\in\R$ such that $z_t<y_t$ and
\begin{equation}\label{eqVeinWag}
\frac{h_t(y_t)-h_t(z_t)}{y_t-z_t}<-\bar{c}.
\end{equation}
In particular, for the inventory problem with incomplete information, in view of convexity of $h$ and the property $h(y-D_t)<+\infty$  for all $y\in\R,$  if the holding cost satisfies condition \eqref{eqVeinWag}, then the functions $h^*_t$ also satisfy this condition.

This is the reason why there is an optimal $(s_t,S_t)$ policy for beliefs.  Of course, this is not an optimal $(s_t,S_t)$ policy for observations.  Based on the history of observations and actions, the DM should compute expected estimations of states, and the optimal $(s_t,S_t)$ policy is based on these estimations.

\bibliography{mainF.bbl}

\end{document}